\documentclass[leqno]{amsart}
\usepackage{amsmath}
\usepackage{amssymb}
\usepackage{amsthm}
\usepackage{enumerate}
\usepackage[mathscr]{eucal}
\theoremstyle{plain}
\usepackage{tikz}
\newtheorem{theorem}{Theorem}[section]

\theoremstyle{definition}
\newtheorem{definition}[theorem]{Definition}
\newtheorem{remark}[theorem]{Remark}

\newtheorem{example}[theorem]{Example}

\theoremstyle{remark}

\begin{document}

\title [Orthogonality and smoothness induced by the norm derivatives] {Orthogonality and smoothness induced by the norm derivatives} 

\author{Debmalya Sain}

\address{(Sain)~Department of Mathematics, Indian Institute of Science, Bengaluru 560012, Karnataka, India}
\email{saindebmalya@gmail.com}

%

\thanks{The author feels elated to acknowledge the delightful friendship of Miss Irin Parvin, while fondly remembering the beautiful times created together.}
\thanks{}

\subjclass[2010]{Primary 46B20,  Secondary 46C50,47L05}
\keywords{norm derivatives; orthogonality in normed spaces; smoothness; linear operators}



\date{}

\begin{abstract}
 We study the concepts of orthogonality and smoothness in normed linear spaces, induced by the derivatives of the norm function. We obtain analytic characterizations of the said orthogonality relations in terms of support functionals in the dual space. We also characterize the related notions of local smoothness and establish its connection with the corresponding orthogonality set, which is analogous to the well-known relation between the Birkhoff-James orthogonality and the classical notion of smoothness. The similarities and the differences between the various notions of smoothness are illustrated by considering some particular examples, including $ \mathbb{K}(\mathbb{H}), $ the Banach space of all compact operators on a Hilbert space $ \mathbb{H}. $ 
\end{abstract}

\maketitle 

\section{Introduction}
The purpose of this article is to study the norm derivatives, from the perspective of two important geometric concepts, namely, orthogonality and smoothness. Let us first establish the notations and the terminologies to be used throughout the article.\\
The letters $\mathbb{X}, \mathbb{Y}, \mathbb{Z}$ stand for normed linear spaces and the letter $ \mathbb{H} $ is used to denote a Hilbert space. The symbol $ \langle~,~\rangle $ is used to denote the inner product on $ \mathbb{H}. $  Let $ \theta $ denote the zero vector of any vector space, other than the scalar field. Throughout the article, we assume that the underlying scalar field is $ \mathbb{R}. $ Let $ B_{\mathbb{X}}= \{ x \in \mathbb{X}~:\|x\| \leq 1 \} $ and $  S_{\mathbb{X}}= \{ x \in \mathbb{X}~:\|x\| = 1 \} $ be the unit ball and the unit sphere of $ \mathbb{X}, $ respectively. Let $\mathbb{L}(\mathbb{X}, \mathbb{Y}) (\mathbb{K}(\mathbb{X}, \mathbb{Y}))$ denote the normed linear space of all bounded (compact) linear operators from $\mathbb{X}$ to $ \mathbb{Y}, $ endowed with the usual operator norm and let $\mathbb X^*$ denote the dual space of $\mathbb{X}.$ We write $ \mathbb{L}(\mathbb{X}, \mathbb{Y}) = \mathbb{L}(\mathbb{X}), $ whenever $ \mathbb{X}=\mathbb{Y}. $ Given $ T \in \mathbb{L}(\mathbb{X}, \mathbb{Y}), $ we use the notations $ \mathcal{R}(T) $ and $ \mathcal{N}(T) $ to denote the range of $ T $ and the kernel of $ T, $ respectively. Let $ M_T = \{ x \in S_{\mathbb{X}} : \| Tx \| = \| T \| \} $ denote the norm attainment set of $ T. $ Given $ x \in \mathbb{X}, $ let $ \mathbb{J}(x) = \{ x^* \in S_{\mathbb{X}^*} : x^*(x) = \|x\| \} $ denote the collection of all support functionals at $ x. $ We say that $ \mathbb{X} $ is smooth at $ x $ if $ \mathbb{J}(x) $ is singleton. It is well-known that smoothness at $ x $ is equivalent to the Gateaux differentiability of the norm at $ x. $ Geometrically, $ \mathbb{X} $ is smooth at $ x $ if and only if there exists a unique supporting hyperplane to the ball $ B(\theta, \|x\|) $ at $ x. $ Since we will also work with some other related notions of smoothness, the above notion of smoothness will sometimes be referred to as the classical notion of smoothness, to avoid any confusion. The study of smoothness in normed linear spaces is intimately connected to the concept of Birkhoff-James orthogonality \cite{B, J}. Given $ x, y \in \mathbb{X}, $ we say that $ x $ is Birkhoff-James orthogonal to $ y, $ written as $ x \perp_B y, $ if $ \| x+\lambda y \| \geq \| x \| $ for all $ \lambda \in \mathbb{R}. $ We use the notation $ x^{\perp} = \{ y \in \mathbb{X} : x \perp_B y \} $ to denote the Birkhoff-James orthogonality set of the vector $ x. $ It is easy to observe that the Birkhoff-James orthogonality is homogeneous, i.e., $ x \perp_B y $ implies that $ \alpha x \perp_B \beta y $ for all  $ \alpha, \beta \in \mathbb{R}. $ We also note that in a Hilbert space, the Birkhoff-James orthogonality relation $ \perp_B $ coincides with the usual orthogonality relation $ \perp $ induced by the underlying inner product $ \langle~,~\rangle. $ As proved in Theorem $ 4.2 $ of the pioneering article by James \cite{J}, a normed linear space $ \mathbb{X} $ is smooth at $ x $ if and only if the Birkhoff-James orthogonality is locally right additive at $ x, $ i.e., $ x \perp_B y, x \perp_B z $ implies that $ x \perp_B (y+z) $ for any $ y, z \in \mathbb{X}. $ \\
The concepts of the norm derivatives arise naturally from the two-sided limiting nature of the Gateaux derivative of the norm, and therefore, they are proper generalizations of the latter. Let us mention the following basic definition of the norm derivatives: 

\begin{definition}
Let $ \mathbb{X} $ be a real normed linear space and let $ x, y \in \mathbb{X}. $ The norm derivatives at $ x $ in the direction of $ y $ are defined as 
\[\rho'_{+} (x, y) = \lim_{t \rightarrow 0+ } \|x\| \frac{\|x+ty\| - \|x\|}{t},\] 

\[\rho'_{-} (x, y) = \lim_{t \rightarrow 0- } \| x \| \frac{\|x+ty\| - \|x\|}{t},\]

\[\rho' (x, y) = \frac{1}{2} \rho{'}_{+} (x, y) + \frac{1}{2} \rho{'}_{-} (x, y) .\]
\end{definition} 

The norm derivatives are extensively applied in studying the geometry of normed linear spaces. We refer the readers to \cite{AST, CW, CWa, M, W, Wa} and the references therein for some of the prominent works in this context. It is worth mentioning that in some cases the function $ \rho' $ may be more convenient than the functions $ \rho'_{+}, \rho'_{-}. $ Indeed, the function $ \rho' $ played a significant role in the paper \cite{Wc}, in order to solve an open problem posed in $ 2010 $ in the book \cite{AST}. For the convenience of the readers, and also for the sake of completeness, we will mention some of the important properties of the norm derivatives in the second section. Our aim in this article is to explore the applicability of the norm derivatives from the geometric perspective of orthogonality and smoothness. Indeed, the above mentioned concepts of $ \rho'_{\pm}, \rho' $ yield natural notions of orthogonality in normed linear spaces \cite{CW, M, Wa}. Given $ x, y \in X, $

\[ x \perp_{\rho_{+}} y \Longleftrightarrow \rho'_{+} (x, y) = 0, \]
 \[ x \perp_{\rho_{-}} y \Longleftrightarrow \rho'_{-} (x, y) = 0, \]
 \[x \perp_{\rho} y \Longleftrightarrow \rho' (x, y) = 0.\]
 
Given $ x \in \mathbb{X}, $ we use the notations $ x^{\perp_{+}}, x^{\perp_{-}}, x^{\perp_{\rho}} $ to denote the $ \rho'_{+} $-orthogonality set, the $ \rho'_{-} $-orthogonality set, and the $ \rho' $-orthogonality set of $ x, $ respectively.\\

We recall that the classical notion of smoothness in a normed linear space is equivalent to the right-additivity of the Birkhoff-James orthogonality. Also, it is not difficult to observe that in general the norm derivatives are not additive in either variable. Moreover, unless the norm is induced by an inner product, there is no hope for the additivity of the norm derivatives in the first variable. Combining all these facts together, the following definition of a generalized smoothness in normed linear spaces turns out to be quite natural. 

\begin{definition}
Let $ \mathbb{X} $ be a real normed linear space and let $ x \in \mathbb{X} \setminus \{\theta\}. $ We say that $ x $ is $ \rho_{+} $-smooth if $ \rho'_{+}(x, .) : \mathbb{X} \longrightarrow \mathbb{R} $ is additive. Similarly, we say that $ x $ is $ \rho_{-} $-smooth, or, $ \rho $-smooth if $ \rho'_{-}(x, .) $ is additive, or,  $ \rho'(x, .) $ is additive, respectively.
\end{definition} 

We would like to note that the concept of $ \rho $-smoothness at a point has been studied previously in \cite{CW, Wa} under the name of semi-smooth points. Moreover, it is not difficult to observe that the above notions of smoothness are proper generalizations of the classical notion of smoothness (the familiar space $ \ell_1 $ serves as an easy example) in a normed linear space.\\

In this article we investigate several notions of orthogonality and smoothness induced by the norm derivatives and their similarities and differences. In particular, we draw parallels between the relation of the classical concept of smoothness with that of the Birkhoff-James orthogonality, and the corresponding counterparts of these concepts induced by the norm derivatives.

\section{Preliminaries}

In this section, mainly for the sake of completeness, we collect a few of the important facts regarding the norm derivatives, some of which find applications in our present work. Since all these results are well-known, we do not provide the proofs, and refer the readers to \cite{A, AST, D, M} for a detailed treatment of the same.\\

\begin{itemize}
\item Using the convexity of the norm function, it is easy to show that the mappings $ \rho'_{\pm} $ (and therefore, the mapping $ \rho' $) are well-defined.\\

\item Given any $ x, y \in \mathbb{X} $ and any $ \alpha \in \mathbb{R}, $ the following holds true:
 \[ (i)~ \rho'_\pm(\alpha x, y)=\rho'_\pm(x,\alpha y) = \left\{
    \begin{array}{ll}
        \alpha \rho'_\pm(x,y),  & \mbox{if } \alpha\geq 0 \\
       \alpha \rho'_\mp(x,y),  & \mbox{if } \alpha< 0.
    \end{array}
\right.\]

In particular, it follows that $ \rho'(\alpha x, y) = \alpha \rho'(x, y) = \rho'(x, \alpha y) $  for all $ \alpha \in \mathbb{R}. $\\

$ (ii)~ \rho'_{-} (x, y) \leq \rho'_{+} (x, y).~ \textit{Moreover, $ \mathbb{X} $ is smooth if and only if}~ \rho'_{+} (x, y) = \rho'_{-} (x, y) $ for all $ x, y \in \mathbb{X}. $ \\

$ (iii)~ \rho'_{\pm} (x, \alpha x + y) = \alpha \| x \|^2 + \rho'_{\pm} (x, y). $ \\

$ (iv)~ \max\{| \rho'_{\pm} (x, y) |, | \rho'(x, y) | \} \leq \| x \| \| y \|. $ \\

$ (v)~ \rho'_{+} (x, y) = \|x\| \sup \{ x^*(y) : x^* \in \mathbb{J}(x) \}. $ \\

$ (vi)~ \rho'_{-} (x, y) = \|x\| \inf \{ x^*(y) : x^* \in \mathbb{J}(x) \}. $ \\

$ (vii)~ $ In particular, given any $ x^* \in \mathbb{J}(x), $ it follows that $ \rho'_{-}(x, y) \leq \| x \| x^*(y) \leq \rho'_{+}(x, y). $ \\

\item The mappings $ \rho'_{\pm}, \rho' $ are continuous with respect to the second variable but not necessarily with respect to the first variable. \\

\item $ x \perp_B y $ if and only if $ \rho'_{-}(x, y) \leq 0 \leq \rho'_{+}(x, y). $ Moreover, either of the conditions $ x \perp_{\rho_{
+}} y $ or $ x \perp_{\rho_{-}} y $ implies that $ x \perp_B y. $
\end{itemize}

\section{Main Results}
 We begin with a characterization of the orthogonality relation $ \perp_{\rho}, $ in terms of supporting functionals.

\begin{theorem}\label{theorem:rho-orthogonality} 

Let $ \mathbb{X} $ be a normed linear space and let $ x, y \in \mathbb{X}.   $ Then $ x \perp_{\rho} y $ if and only if given any $ f \in \mathbb{J}(x), $ there exists $ g \in \mathbb{J}(x) $ such that $ y \in \mathcal{N}(f+g). $

\end{theorem}

\begin{proof}
First we prove the implication $ ``\Longleftarrow". $ We begin with the observation that $ \mathbb{J}(x) $ is a closed subset of $ B_{\mathbb{X}^*} $ with respect to the weak* topology. Therefore, it follows from the Banach-Alaoglou Theorem that $ \mathbb{J}(x) $ is weak*-compact. Combining this with our assumption that given any $ f \in \mathbb{J}(x), $ there exists $ g \in \mathbb{J}(x) $ such that $ y \in \mathcal{N}(f+g), $ we conclude that there exists $ \alpha \geq 0 $ such that $ \{ h(y) : h \in \mathbb{J}(x) \} = [-\alpha, \alpha]. $ From this we deduce the following:
\[ \rho'_{+}(x, y) = \sup_{h \in \mathbb{J}(x)} \|x\| h(y) = \alpha \|x\|, \]
\[ \rho'_{-}(x, y) = \inf_{h \in \mathbb{J}(x)} \|x\| h(y) = -\alpha \|x\|. \] 
Therefore, $ \rho'(x, y) = \frac{1}{2} \rho{'}_{+} (x, y) + \frac{1}{2} \rho{'}_{-} (x, y) = 0, $ as desired.\\

In order to prove $ ``\Longrightarrow" $, let $ f \in \mathbb{J}(x) $ be chosen arbitrarily. Without any loss of generality, we may and do assume that $ f(y) \geq 0. $ Now, we have that $ \rho{'}_{+} (x, y) \geq \|x\| f(y). $ Since $ \rho'_{+} (x, y) + \rho'_{-} (x, y) = 0, $ this implies that $ \rho'_{-} (x, y) \leq -\|x\| f(y). $ As $ \rho'_{-}(x, y) = \inf_{h \in \mathbb{J}(x)} \|x\| h(y), $ there exists $ g_1 \in \mathbb{J}(x) $ such that $ \|x\| g_1(y) \leq -\|x\| f(y) \leq 0, $ which gives that $ g_1(y) \leq -f(y) \leq 0. $ Let $ -f(y) = t g_1 (y), $ where $ t \in [0, 1]. $ If $ t=0, $ then by choosing $ g=f \in \mathbb{J}(x), $ it follows trivially that $ y \in \mathcal{N}(f+g). $ Let us now assume that $ t > 0. $ Let us choose $ m = \frac{2}{1+\frac{1}{t}}. $ Clearly, $ 0 \leq m \leq 1. $ Now, a straightforward computation reveals that 

\begin{align*}
(m g_1 + (1-m) f) (y) & =  \frac{2t}{1+t} g_1 (y) + \frac{1-t}{1+t} f(y) \\
& = - f(y).
\end{align*}

Since $ 1 \geq \| m g_1 + (1-m) f \| \geq | (m g_1 + (1-m) f)x | =1, $ we have essentially proved that $ m g_1 + (1-m) f \in \mathbb{J}(x). $ Choosing $ g = m g_1 + (1-m) f \in \mathbb{J}(x), $ we therefore obtain that $ y \in \mathcal{N}(f+g), $ as desired. This completes the proof of the theorem.

\end{proof}

Recently, in Theorem $ 2.8 $ of \cite{KS}, characterizations of $ \perp_{\rho_{+}} $ and $ \perp_{\rho_{-}} $ have been obtained, in terms of a variation of the Birkhoff-James orthogonality. This, in conjunction with Theorem $ 2.1 $ of \cite{J}, yields the following characterizations of $ \perp_{\rho_{+}} $ and $ \perp_{\rho_{-}}, $ in terms of support functionals. The proof is omitted as it can be completed rather easily by following Theorem $ 2.8 $ of \cite{KS} and then applying Theorem $ 2.1 $ of \cite{J}. It should be noted that the corresponding case of $ \perp_{\rho} $ was not considered in \cite{KS}.

\begin{theorem}\label{theorem:rho(+-)-orthogonality} 

Let $ \mathbb{X} $ be a normed linear space and let $ x, y \in \mathbb{X}.   $ Then the following holds true: \\
(i) $ x \perp_{\rho_+} y $ if and only if there exists $ f_0 \in \mathbb{J}(x) $ such that $ y \in \mathcal{N}(f_0) $ and given any $ f \in \mathbb{J}(x) $ and any $ t \in (0, 1), $ it follows that $ z_t \notin \mathcal{N}(f), $ where $ z_t = -tx+(1-t)y. $\\

(ii) $ x \perp_{\rho_-} y $ if and only if there exists $ f_0 \in \mathbb{J}(x) $ such that $ y \in \mathcal{N}(f_0) $ and given any $ f \in \mathbb{J}(x) $ and any $ t \in (0, 1), $ it follows that $ z_t \notin \mathcal{N}(f), $ where $ z_t = tx+(1-t)y. $

\end{theorem}

Our next result illustrates the relation between the local smoothness induced by the norm derivatives $ \rho'_{\pm} $ and the structure of the corresponding orthogonality sets. As we will note after proving the desired result, it allows us to draw analogies between the classical notion of smoothness and the notion of smoothness induced by the norm derivatives, introduced in the present article.

\begin{theorem}\label{theorem:rho(+-)-smoothness} 

Let $ \mathbb{X} $ be a normed linear space and let $ x \in \mathbb{X}\setminus\{\theta\}. $ Then the following are equivalent: \\
(i) $ x  $ is $ \rho_{+}$-smooth.\\
(ii) $ x  $ is $ \rho_{-}$-smooth.\\
(iii) $ x^{\perp_{+}} $ is a subspace of codimension one in $ \mathbb{X}. $\\
(iv) $ x^{\perp_{-}} $ is a subspace of codimension one in $ \mathbb{X}. $

\end{theorem}

\begin{proof}
Without any loss of generality, we assume that $ \|x\| = 1.$ To prove that $ (i) \implies (ii), $ we begin with arbitrary $ y_1, y_2 \in \mathbb{X}. $ Since $ \mathbb{J}(x) $ is weak*-compact, there exists $ f_0 \in \mathbb{J}(x) $ such that $ \sup_{f \in \mathbb{J}(x)} f(-(y_1+y_2)) = -f_0(y_1+y_2). $ We next observe the following chain of (in)equalities:

\begin{align*}
\rho'_+ (x, -(y_1+y_2)) & = \sup_{f \in \mathbb{J}(x)} f(-(y_1+y_2)) \\
& = f_0(-y_1) + f_0(-y_2) \\
& \leq \rho'_+ (x, -y_1) + \rho'_+ (x, -y_2) \\
& =  \rho'_+ (x, -(y_1 + y_2)),
\end{align*}
where the last equality follows from our assumption that $ x  $ is $ \rho_{+}$-smooth. In particular, this proves that $ \rho'_+ (x, -y_1) = f_0(-y_1) $ and $ \rho'_+ (x, -y_2) = f_0(-y_2). $ Using the properties of  the norm derivatives $ \rho'_{\pm}, $ it is now easy to see that $ \rho'_- (x, y_1) + \rho'_- (x, y_2) = f_0(y_1+y_2). $ Since $ \sup_{f \in \mathbb{J}(x)} f(-(y_1+y_2)) = -f_0(y_1+y_2), $ it follows that $ \inf_{f \in \mathbb{J}(x)} f(y_1+y_2) = f_0(y_1+y_2). $ Therefore, we can deduce that 
\begin{align*}
\rho'_- (x, y_1+y_2) & = \inf_{f \in \mathbb{J}(x)} f(y_1+y_2) \\
& = f_0(y_1 + y_2) \\
& =  \rho'_- (x, y_1)  + \rho'_- (x, y_2).
\end{align*}
Since $ y_1, y_2 \in \mathbb{X} $ were chosen arbitrarily, this proves that $ (i) \implies (ii). $ Using similar arguments, one can now easily prove that $ (ii) \implies (i). $ Let us now prove that $ (i) \implies (iii). $ Given any $ y_1, y_2 \in  x^{\perp_{+}}, $ it follows from $ (i) $ that $ y_1+y_2 \in x^{\perp_{+}}, $ since $ \rho'_+ (x, y_1+y_2) = \rho'_+ (x, y_1) + \rho'_+ (x, y_2) = 0. $ We also observe that
\[ 0=\rho'_+(x, 2y_1-y_1) = 2\rho'_+(x, y_1) - \rho'_-(x, y_1) = - \rho'_-(x, y_1). \]
Now, given any $ \alpha \geq 0, $ clearly, $ \rho'_+ (x, \alpha y_1) = \alpha \rho'_+ (x, y_1) = 0. $ On the other hand, given any $ \alpha < 0, $ it follows from the above observation that $ \rho'_+ (x, \alpha y_1) = \alpha \rho'_-(x, y_1) = 0. $ This proves that $ x^{\perp_{+}} $ is a subspace in $ \mathbb{X}. $ Clearly, $ x^{\perp_{+}} \subsetneq \mathbb{X}, $ as $ x \notin x^{\perp_{+}}. $ Suppose on the contrary to our claim in $ (iii) $ that the codimension of $ x^{\perp_{+}} $ is strictly greater than one.  Then we can find $ y_0 \in \mathbb{X} $ such that $ x, y_0 $ are linearly independent and $ \textit{span}\{x, y_0\} \cap  x^{\perp_{+}} = \{\theta\}. $ Let us choose $ \alpha_0 = - \rho'_+(x, y_0) \neq 0. $ From our assumption that $ x  $ is $ \rho_{+}$-smooth, we obtain by using the properties of $ \rho'_+ $ that
\[ \rho'_+(x, \alpha_0 x + y_0) = \alpha_0 \|x\|^2 + \rho'_+(x, y_0) = 0. \]
Therefore, it follows that  $ \alpha_0 x + y_0 = \theta, $ from our assumption that $ \textit{span}\{x, y_0\} \cap  x^{\perp_{+}} = \{\theta\}. $ However, this clearly contradicts that $ x, y_0 $ are linearly independent. This completes the proof of the assertion that $ (i) \implies (iii). $ We next prove that $ (iii) \implies (i). $ Given any $ y_1, y_2 \in \mathbb{X}, $ there exists $ h_1, h_2 \in x^{\perp_{+}} $ and $ \alpha_1, \alpha_2 \in \mathbb{R} $ such that 
\[ y_1 = \alpha_1 x + h_1 ~\textit{and}~y_2 = \alpha_2 x + h_2. \]
Now, it follows from the properties of $ \rho'_+ $ that
\begin{align*}
\rho'_+ (x, y_1+y_2) & = \rho'_+(x, (\alpha_1+\alpha_2)x + (h_1+h_2)) \\
& = (\alpha_1+\alpha_2) \|x\|^2 + \rho'_+(x, h_1+h_2) \\
& = (\alpha_1+\alpha_2) \|x\|^2 \\
& = \rho'_+(x, y_1) + \rho'_+(x, y_2).
\end{align*}
This completes the proof that $ (iii) \implies (i). $ Finally, we note that the assertion $ (ii) \Longleftrightarrow (iv) $ can be proved by using similar arguments, as given in the above proof of  $ (i) \Longleftrightarrow (iii). $ This establishes the theorem.
\end{proof}

It follows from the above theorem that if $ x \in \mathbb{X} $ is not $ \rho_{\pm} $-smooth then $ x^{\perp_{\pm}} $ is not a subspace of codimension one in $ \mathbb{X}. $ In fact, when $ \mathbb{X} $ is two-dimensional, it is easy to see from Theorem $ 2.8 $ of \cite{KS} that given any $ x \in \mathbb{X} $ which is not $ \rho_{\pm} $-smooth, $ x^{\perp_{+}} $ and $ x^{\perp_{-}} $ are unions of two different sets of rays emanating from the origin. This is further illustrated in the following example.

\begin{example}
Let $ \mathbb{X} = \ell_{\infty}^2 $ and let $ x=(1, 1) \in S_{\mathbb{X}.} $ A simple computation shows that $ x $ is not $ \rho_{\pm} $-smooth and moreover, the following expressions hold true:
\[ x^{\perp_+} = \{(-\alpha, 0) : \alpha \geq 0\} ~\cup~ \{(0, -\beta) : \beta \geq 0 \},\]
\[x^{\perp_-} = \{(\alpha, 0) : \alpha \geq 0\} ~\cup~ \{(0, \beta) : \beta \geq 0 \}.\]
\end{example}

In the same spirit of the above theorem, we mention the following characterization of the local smoothness induced by the norm derivative $ \rho'. $ The proof is omitted, since it can be completed by using the properties of $ \rho' $ and applying similar arguments, as given in the proof of the implication $ (i) \Longleftrightarrow (iii) $ in the above theorem.

\begin{theorem}\label{theorem:rho-smoothness} 

Let $ \mathbb{X} $ be a normed linear space and let $ x \in \mathbb{X}\setminus\{\theta\}. $ Then $ x  $ is $ \rho$-smooth if and only if $ x^{\perp_{\rho}} $ is a subspace of codimension one in $ \mathbb{X}. $
\end{theorem}

In view of the above two theorem, we make the following remarks to point out some important attributes of the various notions of smoothness induced by the norm derivatives.

\begin{remark}
The assertion $ (i) \Longleftrightarrow (ii) $ in Theorem \ref{theorem:rho(+-)-smoothness} shows that the notions of $ \rho_+ $-smoothness and $ \rho_- $-smoothness are identical in any normed linear space. It should be noted that none of them is in general equivalent to the notion of $ \rho $-smoothness. This can be verified readily in case of the space $ \ell_1^2 $ by means of a straightforward computation. We will come back to this in the setting of operator spaces in the later parts of this article. We further observe that it follows from the proof of $ (i) \implies (iii) $ of the same theorem that if $ x \in \mathbb{X} \setminus \{\theta\} $ is $ \rho_+ $-smooth (or, equivalently, $ \rho_- $-smooth) then $ x^{\perp_+} (= x^{\perp_-}) $ is a maximal subspace of $ \mathbb{X}. $
\end{remark}

\begin{remark}
It is well-known (see, for example, \cite{AST, D}) that in a normed linear space $ \mathbb{X}, $ linearity in the second variable of either of the norm derivatives $ \rho'_{\pm} $ is equivalent to the smoothness of $ \mathbb{X}. $ In contrast to such a global result, Theorem \ref{theorem:rho(+-)-smoothness} and Theorem \ref{theorem:rho-smoothness} are essentially local in nature. Of course, the global results follow directly from the corresponding local versions, which further illustrate their applicability.
\end{remark}

\begin{remark}
Theorem \ref{theorem:rho(+-)-smoothness} and Theorem \ref{theorem:rho-smoothness} allow us to draw an interesting analogy between the classical notion of smoothness and the notions of smoothness induced by the norm derivatives $ \rho'_{\pm}, \rho'. $ First of all, it is not difficult to observe that $ x \in \mathbb{X} \setminus \{\theta\} $ is smooth (in the classical sense) if and only if $ x^{\perp} = \{ y \in \mathbb{X} : x \perp_B y \} $ is a maximal subspace of $ \mathbb{X}. $ Although the notions of $ \rho_{\pm} $-smoothness and $ \rho $-smoothness are not equivalent to the classical notion of smoothness, this particular characterization remains valid in each of these cases. Indeed, as the above two theorems state, $ x \in \mathbb{X} \setminus \{\theta\} $ is $ \rho_{\pm} $-smooth ($ \rho $-smooth) if and only if $ x^{\perp_\pm} (x^{\perp_\rho}) $ is a maximal subspace of $ \mathbb{X}. $
\end{remark}

We would like to mention that for the remaining part of this article, which deals with $ \rho_{\pm} $-smoothness in operator spaces, Theorem $ 3.2 $ and Theorem $ 3.4 $ of \cite{W} play a central role. We begin by characterizing the $ \rho_\pm $-smooth points of the space of all compact operators on a Hilbert space, endowed with the usual operator norm.

\begin{theorem}\label{theorem:rho(+-)-smoothness in K(H)} 

Let $ \mathbb{H} $ be a Hilbert space and let $ T \in \mathbb{K}(\mathbb{H}) $ be non-zero. Then the following are equivalent: \\
(i) $ T  $ is $ \rho_+$-smooth. \\
(ii) $ T  $ is $ \rho_-$-smooth. \\
(iii) $ T $ is smooth. \\
(iv) $ M_T = \{ \pm x_0 \}, $ for some $ x_0 \in S_{\mathbb{H}}. $
\end{theorem}

\begin{proof}
Without any loss of generality, we assume that $ \|T\| = 1. $ We begin with the following three basic observations:\\
$ (a) $ The equivalence of $ (i) $ and $ (ii) $ follows from Theorem \ref{theorem:rho(+-)-smoothness}.\\
$ (b) $ The equivalence of $ (iii) $ and $ (iv) $ follows from Theorem $ 4.1 $ and Theorem $ 4.2 $ of \cite{PSG}. \\
$ (c) $ The implication $ (iii) \implies (i) $ follows trivially.\\

Therefore, to complete the proof of the theorem, it is sufficient to show that $ (i) \implies (iv). $ It follows from Theorem $ 2.2 $ of \cite{SP} and the fact that $ T $ is compact, that $ M_T = S_{\mathbb{H}_0}, $ for some finite-dimensional subspace $ \mathbb{H}_0 $ of $ \mathbb{H}. $ If $ \mathbb{H}_0 $ is one-dimensional then we are done. Suppose on the contrary that $ \dim \mathbb{H}_0 > 1. $ Let $ \{ e_1, \ldots, e_n \} $ be an orthonormal basis of $ \mathbb{H}_0, $ where $ n > 1. $ We extend it to a complete orthonormal basis $ \mathcal{B} = \{ e_{\alpha} : \alpha \in \Lambda \} $ of $ \mathbb{H}, $ where $ \Lambda $ is the index set. Let us now define $ A_1, A_2 \in \mathbb{K}(\mathbb{H}) $ in the following way:

\begin{align*}
\noindent A_1(e_1) = Te_1 & \hspace{3cm} A_2(e_1) = -\frac{1}{2}Te_1 \\
\noindent A_1(e_2) = -\frac{1}{2}Te_2 & \hspace{3cm} A_2(e_2) = Te_2 \\
\noindent A_1(e_\alpha) = 0 & \hspace{3cm} A_2(e_\alpha) = 0 ~\forall~ \alpha \in \Lambda \setminus \{ 1, 2 \}. 
\end{align*}

Applying Theorem $ 3.4 $ of \cite{W}, we obtain the following:\\
\[ \rho'_+(T, A_1) = \max\{\langle Tx, A_1x \rangle : x \in M_T\} = \|T\|^2 = 1, \]
\[ \rho'_+(T, A_2) = \max\{\langle Tx, A_2x \rangle : x \in M_T\} = \|T\|^2 =1. \]
On the other hand, given any $ z = \sum_{i=1}^{n} \alpha_i e_i \in M_T, $ a trivial computation shows that 
\[ \langle Tz, (A_1+A_2)z \rangle =  \frac{1}{2} (\alpha_1^2+\alpha_2^2 ) \|T\|^2 \leq \frac{1}{2}, \]
which shows that $ \rho'_+(T, A_1+A_2) \neq \rho'_+(T, A_1) + \rho'_+(T, A_2), $ contradicting our assumption that $ T  $ is $ \rho_+$-smooth. Therefore, it must be true that $ \mathbb{H}_0 $ is one-dimensional. This completes the proof of the theorem. 

\end{proof}

In light of the above theorem, we would like to make the following remark which illustrates that the concepts of $ \rho_{\pm} $-smoothness are not equivalent to that of $ \rho $-smoothness in a general normed linear space. We also note that it follows from the above theorem that in $ \mathbb{K}(\mathbb{H}), $ $ \rho_{\pm} $-smoothness is equivalent to the classical smoothness.

\begin{remark}
Recently it has been proved in \cite{Wb} that given any real Hilbert space $ \mathbb{H}, $ an operator $ T \in \mathbb{K}(\mathbb{H}) $ is $ \rho $-smooth if and only if $ \dim (\textit{span}~M_T) \leq 2. $ We would like to note that in case of the two-dimensional Euclidean plane, this can be verified directly by applying Theorem $ 3.4 $ of \cite{W} and a straightforward computation. In particular, it now follows from Theorem \ref{theorem:rho(+-)-smoothness in K(H)} and the above observation that $ \rho_{\pm} $-smoothness and $ \rho $-smoothness are not equivalent in $ \mathbb{K}(\mathbb{H}). $
\end{remark}

In general, it seems to be a difficult problem to explicitly identify the $ \rho_{\pm} $-smooth points of $ \mathbb{L}(\mathbb{X}, \mathbb{Y}), $ where $ \mathbb{X}, \mathbb{Y} $ are given normed linear spaces. However, it is easy to see that by applying Theorem $ 3.2 $ of \cite{W}, the following sufficient condition for $ \rho_{\pm} $-smoothness can be obtained in the setting of $ \mathbb{L}(\mathbb{X}, \mathbb{Y}), $ under some additional assumptions. We refer the readers to \cite{HWW} for an excellent exposition on the theory and applications of M-ideals in Banach spaces.

\begin{theorem}\label{theorem:rho(+-)-smoothness in K(X, Y)} 
Let $ \mathbb{X} $ be a reflexive Banach space and let $ \mathbb{Y} $ be any Banach space such that $ \mathbb{K}(\mathbb{X}, \mathbb{Y}) $ is an M-ideal in $ \mathbb{L}(\mathbb{X}, \mathbb{Y}). $ Let $ T \in \mathbb{L}(\mathbb{X}, \mathbb{Y}) $ be such that dist$ (T, \mathbb{K}(\mathbb{X}, \mathbb{Y})) < \| T \| $ and $ M_T = \{ \pm x_0 \}, $ for some $ x_0 \in S_{\mathbb{X}}. $ Let $ Tx_0 $ be a $ \rho_{\pm}$-smooth point in $ \mathbb{Y}. $  Then $ T $ is a $ \rho_{\pm} $-smooth point in $ \mathbb{L}(\mathbb{X}, \mathbb{Y}). $
\end{theorem}

\begin{proof}
Without any loss of generality, we assume that $ \| T \|=1. $ Let $ A_1, A_2 \in \mathbb{L}(\mathbb{X}, \mathbb{Y}) $ be arbitrary. Since $ M_T = \{ \pm x_0 \}, $ Theorem $ 3.2 $ of \cite{W} asserts that 
\[\rho'_+(T, A_1) = \rho'_+(Tx_0, A_1x_0),~~ \rho'_+(T, A_2) = \rho'_+(Tx_0, A_2x_0).\]
On the other hand, using the $ \rho_+ $-smoothness of $ Tx_0, $ we obtain by applying the same theorem and the above expressions that
\[ \rho'_+(T, A_1+A_2) = \rho'_+(Tx_0, (A_1+A_2)x_0) = \rho'_+(T, A_1)  +\rho'_+(T, A_2).\]

This proves that $ T $ is $ \rho_+ $-smooth. Combining this with Theorem \ref{theorem:rho(+-)-smoothness}, we obtain the desired conclusion that $ T $ is $ \rho_{\pm} $-smooth. This completes the proof of the theorem.
\end{proof}

Using the expression of the norm derivative $ \rho', $ it is straightforward to observe that the above theorem also provides a sufficient condition for $ \rho $-smoothness, under the modified assumption that $ Tx_0 $ is a $ \rho $-smooth point in $ \mathbb{Y}. $ In particular, it is easy to conclude that $ \rho_{\pm} $-smoothness and $ \rho $-smoothness are not equivalent in $ \mathbb{L}(\mathbb{X}, \mathbb{Y}). $ This is recorded in the following theorem, and we omit the proof as it is now clear in view of the above comment.

\begin{theorem}\label{theorem:rho-smoothness in K(X, Y)}
Let $ \mathbb{X} $ be a reflexive Banach space and let $ \mathbb{Y} $ be any Banach space such that $ \mathbb{K}(\mathbb{X}, \mathbb{Y}) $ is an M-ideal in $ \mathbb{L}(\mathbb{X}, \mathbb{Y}). $ Let $ T \in \mathbb{L}(\mathbb{X}, \mathbb{Y}) $ be such that dist$ (T, \mathbb{K}(\mathbb{X}, \mathbb{Y})) < \| T \| $ and $ M_T = \{ \pm x_0 \}, $ for some $ x_0 \in S_{\mathbb{X}}. $ Let $ Tx_0 $ be a $ \rho$-smooth point in $ \mathbb{Y}. $ Then $ T $ is a $ \rho $-smooth point in $ \mathbb{L}(\mathbb{X}, \mathbb{Y}). $
\end{theorem}

Let us end this article with the following complete characterization of $ \rho_{\pm} $-smooth points in $ \mathbb{L}(\ell_{1}^n, \mathbb{Y}), $ for any normed linear space $ \mathbb{Y}. $

\begin{theorem}\label{theorem:L(l_1^n, Y)} 
Let $ \mathbb{Y} $ be any normed linear space and let $ T \in \mathbb{L}(\ell_{1}^n, \mathbb{Y}). $ Then $ T $ is $ \rho_{\pm} $-smooth if and only if $ M_{T} = \{ \pm u_0 \}, $ for some extreme point $ u_0 $ of $ B_{\ell_{1}^n} $ and $ Tu_0 $ is a $ \rho_{\pm} $-smooth point of $ \mathbb{Y}. $
\end{theorem}

\begin{proof}
Since the implication $ ``\Longleftarrow" $ follows directly from Theorem \ref{theorem:rho(+-)-smoothness in K(X, Y)}, let us prove only the implication $ ``\Longrightarrow". $ Let $ \{ e_i : 1 \leq i \leq n \} $ denote the usual canonical basis of $ \ell_1^n. $ Clearly, the set of extreme points of $ B_{\ell_{1}^n} $ are given by $ \{ \pm e_i : 1 \leq i \leq n \}. $ Suppose on the contrary that there exist $ e_i, e_j \in M_T $ such that $ i \neq j. $ Without any loss of generality, let us assume that $ i=1 $ and $ j=2. $ Define $ A_1, A_2 \in \mathbb{L}(\ell_{1}^n, \mathbb{Y}) $ in the following way:

\begin{align*}
\noindent A_1(e_1) = Te_1 & \hspace{3cm} A_2(e_1) = -\frac{1}{2}Te_1 \\
\noindent A_1(e_2) = -\frac{1}{2}Te_2 & \hspace{3cm} A_2(e_2) = Te_2 \\
\noindent A_1(e_m) = 0 & \hspace{3cm} A_2(e_m) = 0 ~\forall~ m \in \{3, \ldots, n\}. 
\end{align*}

As in the proof of Theorem \ref{theorem:rho(+-)-smoothness in K(H)}, it is easy to check that $ \rho'_+(T, A_1+A_2) \neq \rho'_+(T, A_1) + \rho'_+(T, A_2), $ contradicting our assumption that $ T  $ is $ \rho_+$-smooth. This proves that $ M_{T} = \{ \pm u_0 \}, $ for some extreme point $ u_0 $ of $ B_{\ell_{1}^n}. $ We next claim that $ Tu_0 $ is a $ \rho_{\pm} $-smooth point of $ \mathbb{Y}. $ If this claim is not true then there exists $ y_1, y_2 \in \mathbb{Y} $ such that $ \rho'_+(Tu_0, y_1+y_2) \neq \rho'_+(Tu_0, y_1) + \rho'_+(Tu_0, y_2). $ Let $ A_3, A_4 \in \mathbb{L}(\ell_{1}^n, \mathbb{Y}) $ be such that $ A_3u_0=y_1 $ and $ A_4u_0 = y_2. $ It is easy to check by applying Theorem $ 3.2 $ of \cite{W} that $ \rho'_+(T, A_3+A_4) \neq \rho'_+(T, A_3) + \rho'_+(T, A_4), $ once again contradicting our assumption that $ T  $ is $ \rho_+$-smooth. This prove that our claim of $ Tu_0 $ being a $ \rho_{\pm} $-smooth point of $ \mathbb{Y} $ is correct. Hence the theorem.
\end{proof}

\bibliographystyle{amsplain}

\end{document}